\documentclass[reqno]{amsart}

\usepackage{amsfonts,amssymb,amsmath}
\usepackage{color}
\usepackage{graphicx}

\newtheorem{theorem}{Theorem}[section]
\newtheorem{lemma}[theorem]{Lemma}
\newtheorem{proposition}[theorem]{Proposition}
\newtheorem{corollary}[theorem]{Corollary}
\newtheorem{remark}[theorem]{Remark}
\newtheorem{definition}{Definition}[section]

\newcommand{\mc}[1]{{\mathcal #1}}
\newcommand{\mf}[1]{{\mathfrak #1}}

\newcommand{\bb}[1]{{\mathbb #1}}

\newcommand{\eps}{\varepsilon}

\newcommand{\<}{\langle}
\renewcommand{\>}{\rangle}
\newcommand{\p}{\partial}
\newcommand{\pfrac}[2]{\genfrac{}{}{}{1}{#1}{#2}}

\newcommand{\dl}{\<\!\<}
\newcommand{\dr}{\>\!\>}

\title[Hydrodynamic Limit]{Hydrodynamic Limit for a Type of
Exclusion Processes with slow bonds in dimension $\ge 2$}

\author[T. Franco]{Tertuliano Franco}
\address{Tertuliano Franco \hfill\break\indent 
IMPA \hfill\break\indent 
Estrada Dona Castorina 110, \hfill\break\indent
J. Botanico, 22460 Rio
de Janeiro, Brazil} 
\email{tertu@impa.br}

\author[A. Neumann]{Adriana Neumann}
\address{Adriana Neumann \hfill\break\indent 
IMPA \hfill\break\indent 
Estrada Dona Castorina 110, \hfill\break\indent
J. Botanico, 22460 Rio
de Janeiro, Brazil} 
\email{aneumann@impa.br}

\author[G. Valle]{Glauco Valle}
\address{Glauco Valle \hfill\break\indent 
UFRJ \hfill\break\indent 
Departamento de M\'etodos Estat\'isticos do Instituto de Matem\'atica, \hfill\break\indent
Caixa Postal 68530, 21495-970, Rio de Janeiro, Brazil} 
\email{glauco.valle@im.ufrj.br}

\date{\today}

\begin{document}

\maketitle

\begin{abstract} Let $\Lambda$ be a connected closed region with smooth boundary contained in the
$d$-dimensional continuous torus $\bb T^d$. In the discrete torus $N^{-1} \bb T^d_N$, we consider a nearest neighbor symmetric
exclusion process where occupancies of neighboring sites are exchanged at rates depending on $\Lambda$ 
in the following way: if both sites are in $\Lambda$ or $\Lambda^\complement$, 
the exchange rate is one; If one site is in $\Lambda$ and the other one is in $\Lambda^\complement$ 
and the direction of the bond connecting the sites is $e_j$, then the exchange rate
is defined as $N^{-1}$ times the absolute value of the inner product between $e_j$ and the normal exterior vector to
$\p\Lambda$. We show that this exclusion type process has a non-trivial 
hydrodynamical behavior under diffusive scaling and, in the continuum limit, particles
are not blocked or reflected by $\partial\Lambda$. Thus the model represents a system of particles
under hard core interaction in the presence of a permeable membrane which slows down the passage 
of particles between two complementar regions. 
\end{abstract}

\section{Introduction}

The exclusion process is a continuous time interacting particle system where particles 
move as independent random walks on a graph except for the exclusion rule that prevents 
two particles from occupying the same site, or vertex. In the symmetric case, the 
process evolves as follows: to each bond we associate a waiting exponential time,
which are independent of the waiting time for any other bond; at the waiting time the 
occupancies of the sites connected by the bond are exchanged; the parameter of the
exchange times, or exchange rate, depends only on the bond. The especification of
the exchange rates determines the environments for the exclusion process. 
In our case, as the underlying graph, we consider the discrete torus with $N^d$ points
and nearest neighbor bonds. The variable $N$ is the scaling parameter.

This paper studies the hydrodynamical behavior of symmetric exclusion processes in 
non-homogeneous environments with the presence of slow bonds. Here we mean non-constant 
environments where a usual bond has exchange rate one and a slow bond has 
exchange rate lower than one. With respect to the scaling parameter, we assume that a slow bond
has exchange rate of order $N^{-1}$. When the environment is constant, the exclusion process
has a well-known hydrodynamical behavior under diffusive scaling, but, in the presence of slow 
bonds, particles will not move fast enough to garantee that we still have an hydrodynamic 
behavior in diffusive scaling. Hydrodynamics in diffusive scaling have been obtained in
several cases, even when the environment is random and consists only of slow bonds.    

For one dimensional processes, in \cite{fjl}, the exchange rate over a bond $[\frac{x}{N},\frac{x+1}{N}]$ is 
given by $[N(W(x+1/N)-W(x/N))]^{-1}$, where $W$ is an $\alpha$-stable subordinator of a 
L\'evy Process. They obtain a quenched hydrodynamic limit. In papers previous 
to \cite{fjl}, for example \cite{f} and \cite{n}, the randomness or non-homogenity did not survive 
in the continuum limit. An also one-dimensional result, following \cite{fjl}, was 
obtained in \cite{fl}, for more general, but non-random, increasing functions $W$. 
The techniques used in those papers were strongly based on theorems about convergence of one dimensional 
continuous time stochastic processes. In fact, even the $d$-dimensional case treated in 
\cite{v} has considered a class of non-homogeneous environments that could be decomposed, in a 
proper sense, in $d$ one-dimensional cases. Recently, different approaches have been searched 
to deal with $d$-dimensional environments, see \cite{f2} and \cite{j}.

We now describe the exclusion processes we are concerned. Let $\{e_j:j=1,...,d\}$ be the 
canonical base of $\mathbb{R}^d$ and $\Lambda\subset\bb T^d$ be a simple connected region 
with smooth boundary $\partial\Lambda$. If the bond $[\frac{x}{N},\frac{x+e_j}{N}] \in N^{-1} \bb T^d_N$ has vertices in 
each of the regions $\Lambda$ and $\Lambda^\complement$, its exchange rate is defined as 
$N^{-1}$ times the absolute value of the inner product between $e_j$ and the normal exterior vector to
$\p\Lambda$. For others edges, the exchange rate is defined as one. This means that the 
slow bonds are among those crossing the boundary of $\Lambda$. We call this process the exclusion process
with slow bonds over $\p\Lambda$.

We can interpret $\p\Lambda$ as a permeable membrane, which slows down the passage of particles between 
the regions $\Lambda$ and $\Lambda^\complement$. For this type of exclusion process, 
the membrane does not completely prevent the passage of particles, and still survives in 
the continuum limit, appearing explicitely in the hydrodynamic equation. The exchange rate of particles 
for a bond crossing $\p\Lambda$ is smaller if the bond is close to a tangent line of $\p\Lambda$. 
Note that this assumption has physical meaning, take for example cases of reflections in several 
physical models: partial reflection of light crossing a media with diferent refraction indexes, 
mechanical systems where particles try to cross some interface, etc. However the direction of 
the speed of particles is not changed as usually occur in physical reflection. Our 
definition of the exchange rates also allows a strong convergence result for the 
empirical measures associated to the exclusion process making simpler the 
proof of the hydrodynamic limit. 

The hydrodynamical equation of the exclusion process with slow bonds over $\p\Lambda$ is a 
parabolic partial differential equation $\displaystyle \partial_t \rho \; =\; \mc L_\Lambda \rho$, where 
the operator $\mc L_\Lambda$ is a sort of $d$-dimensional Krein-Feller operator. Without the presence 
of slow bonds, the operator $\mc L_\Lambda$ would be replaced by the laplacian operator acting on $C^2$ 
functions and the hydrodynamical equation is therefore the heat equation. Here, the existence of the 
membrane modifies the domain, and thus the operator itself. In fact, we observe that the proper domain for 
$\mc L_\Lambda$ contains functions that are discontinuous over $\p\Lambda$. Geometrically, $\mc L_\Lambda$ 
glues the discontinuity of a function around $\p\Lambda$ and then behaves like the laplacian.

One possible approach to prove the hydrodynamic limit for the exclusion process with slow bonds over 
$\p\Lambda$ is through Gamma convergence. In \cite{j}, this approach and the conditions for it to hold
are discussed, see also \cite{f}. There, the coersiveness condition would require some kind of Rellich-Kondrachov's Theorem 
(namely, the compact embedding in $L^2$ of some sort of Sobolev space supporting an extension of $\mc L_\Lambda$, 
see \cite{e}). In the method presented here, we go in this direction, but instead of reach the hypotheses 
in \cite{j}, we have use similar the analytical tools to obtain a short and simple proof of uniqueness of the
hydrodynamic equation. We also show that the extension of $\mc L_\Lambda$ satisfies the Hille-Yoshida 
Theorem. On the other hand, the convergence from discrete to continuous that we present here is made in very 
direct way, and it was inspired by the convergence of the discrete laplacian to the continuous laplacian. 

\smallskip

The paper is presented as follows:
In Section \ref{2}, we define the model and state all results proved to be proved in the paper;
Section \ref{3} is devoted to prove all the results concerning to the continuous operator $\mc L_\Lambda$;
In Section \ref{4}, the hydrodynamic limit is proved.

\section{Notation and Results}\label{2}

Let $\bb T^d$ be the $d$-dimensional torus, which is $[0,1)^d$ with periodic boundary conditions,
and $\bb T^d_N$ be the discrete torus with $N^d$ points,  which is $\{0,...,N-1\}^d$ with periodic 
boundary conditions. We denote by $\eta = (\eta(x))_{x \in \bb T_{N}^d}$ a typical configuration 
in the state space $\Omega_N = \{0,1\}^{\bb T_{N}^d}$, for which, $\eta(x)=0$ means that
site $x$ is vacant, and $\eta(x)=1$ that site $x$ is occupied. If a bond of $N^{-1} \bb T^d_N$ has
vertices $\pfrac{x}{N}$ and $\pfrac{y}{N}$, it will be denoted by $[\frac{x}{N},\frac{y}{N}]$.

Recall that $\{e_j:j=1,...,d\}$ is the canonical base of $\mathbb{R}^d$. The symmetric nearest neighbor 
exclusion process with exchange rates $\xi^{N}_{x,y} > 0$, $x,y \in \bb T^d_N$, $|x-y| = 1$, is a 
Markov Process with configuration space $\Omega_N$, whose generator $L_{N}$ acts on functions 
$f:\Omega_N \rightarrow \bb{R}$ as
\begin{equation}\label{ln}
(L_{N}f)(\eta)=\sum_{x\in \bb T_{N}^d}\sum_{j=1}^d \,\xi^{N}_{x,x+e_j}\,\Big[f(\eta^{x,x+e_j})-f(\eta)\Big]\,,
\end{equation}
where $\eta^{x,x+e_j}$ is the configuration obtained from $\eta$ by exchanging the variables 
$\eta(x)$ and $\eta(x+e_j)$:
\begin{equation*}
(\eta^{x,x+e_j})(y)=\left\{\begin{array}{cl}
\eta(x+e_j),& \mbox{if}\,\,\, y=x\,,\\ 
\eta(x),& \mbox{if} \,\,\,y=x+e_j\,,\\ 
\eta(y),& \mbox{otherwise.}
\end{array}
\right.
\end{equation*}
Let $\nu^N_\alpha$, $\alpha \in (0,1)$, be the Bernoulli product measure $\Omega_N$, i.e., the product measure
whose marginals have Bernoulli distribution with parameter $\alpha$. Then $\{\nu^N_\alpha : 0\le \alpha \le 1\}$ 
is a family of invariant, in fact reversible, measures for any symmetric exclusion process.
\smallskip

Now, fix a simple connected region $\Lambda\subset\bb T^d$ with smooth boundary $\partial\Lambda$.
Denote by $\vec{\zeta}(u)$ the normal unitary exterior vector to the smooth surface $\p\Lambda$ in the
point $u\in\p\Lambda$. If $\pfrac{x}{N}\in\Lambda$ and $\pfrac{x+e_j}{N}\in\Lambda^\complement$, or $\pfrac{x}{N}\in\Lambda^\complement$ and $\pfrac{x+e_j}{N}\in\Lambda$, we
define $\vec{\zeta}_{x,j}$ as a vector $\vec{\zeta}(u)$ evaluated in an arbitrary but fixed point 
$u\in \p\Lambda \cap [x,x+e_j]$. The exclusion process with slow bonds over $\p\Lambda$ is a
symmetric nearest neighbor exclusion process with exchange rates $\xi^N_{x,x+e_j}\;=\;\xi^N_{x+e_j,x}$ given by
\begin{equation}
\label{rates}
\left\{\begin{array}{cl}
\dfrac{\vert \vec{\zeta}_{x,j}\cdot e_j\vert}{N}\,, &  \mbox{if}\,\,\,\,\frac{x}{N}\in\Lambda\textrm{~and~}
\frac{x+e_j}{N}\in\Lambda^\complement\,, \textrm{ or } \pfrac{x}{N}\in\Lambda^\complement \textrm{ and } \pfrac{x+e_j}{N}\in\Lambda\,,\\ 
\quad\\
1\,, &\mbox{otherwise,}
\end{array}
\right.
\end{equation}
for $j=1,\ldots,d$, and for every $x \in \bb T^d_N$. In this case, the exchange rate of a bond crossing 
the boundary $\p\Lambda$ is also of order $N^{-1}$, but it depends on the angle of incidence: the crossing
of $\p\Lambda$ by a particle gets harder to happen as the direction of entrance gets closer to the tangent 
plane to the surface $\p\Lambda$.

\begin{figure}[t]
 \centering
\input{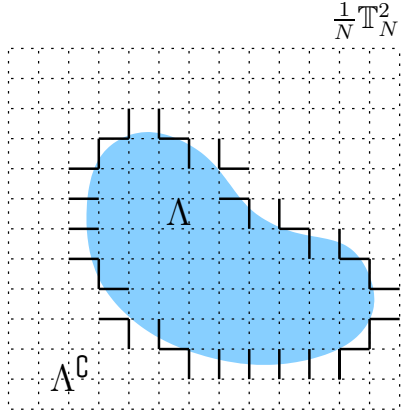}
 \caption{ The darker region corresponds to $\Lambda$. The bolded bonds have exchanges rates 
$\frac{\vert\vec{\zeta}_{x,j}\cdot e_j \vert}{N}$, any other bond has exchange rate $1$.}
 \label{fig5}
\end{figure}

\smallskip

From now on, the rates in the definition of $L_N$ will always be given by \ref{rates}.
Denote by $\{\eta^N_t : t\ge 0\}$ a Markov process with state space $\Omega_N$ and generator 
$L_N$ speeded up by $N^2$. Let $D(\bb R_+, \Omega_N)$ be the Skorohod space of
c\`adl\`ag trajectories taking values in $\Omega_N$. For a measure $\mu$ on $\Omega_N$, denote 
by $\bb P^N_{\mu}$ the probability measure on $D(\bb R_+, \Omega_N)$ induced by the
initial state $\mu$ and the Markov process $\{\eta_t^N : t\ge 0\}$. The expectation with 
respect to $\bb P^N_{\mu}$ is going to be denoted by $\bb E^N_{\mu}$.

\smallskip

A sequence of probability measures $\{\mu_N : N\geq 1 \}$ is said to be associated to a profile 
$\gamma :\bb T^d \to [0,1]$ if $\mu_N$ is a probability measure on $\Omega_N$, for every N, and
\begin{equation}
\label{f09}
\lim_{N\to\infty}
\mu_N \left\{ \, \Big\vert \pfrac {1}{N^d}\!\! \sum_{x\in\bb T_N^2} H(\pfrac{x}{N}) \eta(x)
- \int H(u) \gamma(u) du \Big\vert > \delta \right\} \;=\; 0
\end{equation}
for every $\delta >0$, and every continuous functions $H: \bb T^d \to \bb R$.

\smallskip

The exclusion process with slow bonds over $\p\Lambda$ has a related random walk on $N^{-1} \bb T^d_N$ 
that describes the evolution of the system with a single particle. Thus particles in the exclusion 
process evolve independently as such random walk except for the hard core interaction. To simplify 
notation later, we introduce here the generator of this random walk, which is given by
\begin{equation*}
(\mathbb{L}_N H)(\pfrac{x}{N}) = \sum_{j=1}^d \Big\{ \xi^N_{x,x+e_j} \, \Big[ H(\pfrac{x+e_j}{N}) 
- H(\pfrac{x}{N}) \Big] + \xi^N_{x,x-e_j} \, \Big[H(\pfrac{x-e_j}{N}) - H(\pfrac{x}{N}) \Big] \Big\}\,,
\end{equation*}
for every $H: N^{-1} \bb T^d_N \rightarrow \mathbb{R}$ and every $x \in \bb T^d_N$. We will not 
differentiate the notation for functions $H$ defined on $\bb{T}^d$ and on $ N^{-1} \bb T_{N}^d$.


\medskip

\subsection{The Operator $\mc L_\Lambda$}

Here we define the operator $\mc L_\Lambda$ and state its main properties. Its domain is 
defined as a set of functions two times continuously differentiable inside and outside $\Lambda$, 
but possibly discontinuous on $\partial \Lambda$. Besides, in the boundary $\partial \Lambda$, 
we should require particular conditions of those functions in order to have good properties 
of $\mc L_\Lambda$ that allows us to conclude the uniqueness of solutions of the hydrodynamic 
equation and obtain a strong convergence result for the empirical measures in the proof of 
the hydrodynamic limit. The necessity of these conditions are going to be made clear later in the text.

\smallskip

\begin{definition}\label{operator} Recall that $\vec{\zeta}$ denotes the normal exterior vector to the
surface $\p\Lambda$. The domain $\mf D_\Lambda\subset L ^2(\bb T^d)$ will be the set of functions 
 $H\in L^2(\bb T^d)$, such that $H(u)=h(u)+\lambda {\bf 1}_\Lambda(u)$, 
where:
\begin{itemize}
\item[\textbf{(i)}] $\lambda\in \bb R$;
\item[\textbf{(ii)}] $h\in C^2(\bb T^d)$;
\item[\textbf{(iii)}] $\nabla h \vert_{\p\Lambda}(u)=-\lambda\,\vec{\zeta}(u)$.

\end{itemize}
Now, we define the operator  $\mc L_\Lambda: \mf D_\Lambda\rightarrow L^{2}(\bb T^d)$ by
\begin{equation*}
 \mc L_\Lambda H\,=\Delta h\,.
\end{equation*}
\end{definition}

Geometrically, the operator $\mc L_\Lambda$ removes the discontinuity around the
surface $\partial\Lambda$ and then acts like the laplacian operator.

\begin{remark}
It is not entirely obvious why there exist functions $h\in C^2(\bb T^d)$ such that
$\nabla h \vert_{\p\Lambda}(u)=-\lambda\,\vec{\zeta}(u)$, for $\lambda \neq 0$. For an example of such a function,
consider firstly  $g:\bb T^d\to \bb R$ defined by
\begin{equation*}
g(u)\;=\;\left\{\begin{array}{cl}
\lambda \, \textrm{dist}\,(u,\p\Lambda)\,, &  \mbox{if}\,\,\,\,u\in\Lambda^\complement\,,\\
- \lambda \, \textrm{dist}\,(u,\p\Lambda)\,, &  \mbox{if}\,\,\,\,u\in\Lambda\,.
\end{array}
\right.
\end{equation*}
One can checks that this function is smooth in an open neighborhood $V$ of $\p\Lambda$, and
satisfies the condition $\nabla g \vert_{\p\Lambda}(u)=-\lambda\,\vec{\zeta}(u)$. 
However, $g$ is  not differentiable in the space $\bb T^d$. To solve this problem,
it is enough to multiply $g$ by $\sum_i \Phi_i$, where $\{\Phi_i\}$ is a partition of unity
such that the support of any  $\Phi_i$ is contained in $V$ and $\sum_i \Phi_i(u)=1$ for
all $u\in U\subset V$, $U$ an open set containing $\p\Lambda$. Finally, the function
\begin{equation*}
 h(u)\,=\,g(u)\,\sum_i \Phi_i(u)
\end{equation*}
satisfies the required conditions.
\end{remark}

For the next result we need to introduce some notation. We denote by $\bb I$ the identity
operator in $L^2(\bb T^d)$ and by $\dl \cdot , \cdot \dr$ and $\| \cdot \|$ its usual inner product and norm: 
\begin{equation*}
\dl f, g \dr \;=\; \int_{\bb T^d} f(u)\, g(u)\, du\; \, \textrm{ and } \, 
\| f \| = \sqrt{ \dl f, f \dr } \; , \ f, \ g \in L^2(\bb T^d)\,.
\end{equation*}

\begin{theorem}
\label{Lambda}
There exists a Hilbert Space $(\mc H^1_\Lambda,\dl \cdot, \cdot \dr_{1,\Lambda})$ compactly embedded in $L^2(\bb T^d)$ such that $\mf D_\Lambda \subset \mc H^1_\Lambda$ and $\mc L_\Lambda$ can be extended to $\mc L_\Lambda: \mc H^1_\Lambda  \rightarrow L^2(\bb T^d)$ in such way that the extension enjoys the following properties:
\begin{itemize}
\item[(a)] The domain $\mc H^1_\Lambda $ is dense in $L^2(\bb T^d)$;
\item[(b)] The operator $\mc L_\Lambda$ is self-adjoint and non-positive:
$\dl H, -\mc L_\Lambda H\dr \geq 0$, for all $H$ in $\mc H^1_\Lambda $;
\item[(c)] The operator $\bb I -\,\mc L_\Lambda : \mc H^1_\Lambda \rightarrow L^2(\bb T^d)$ is bijective and $\mf D_\Lambda$ is a core for it;
\item[(d)] The operator $\mc L_\Lambda$ is dissipative, i.e., 
$$
\Vert \mu H -\mc L_\Lambda H\Vert \geq \mu \Vert H\Vert \, ,
$$
for all $H\in\mc H^1_\Lambda $ and $\mu >0$;
\item[(e)] The eigenvalues of $- \mc L_\Lambda$ form a countable set
$0= \mu_0\le \mu_1\le \cdots$ with $\lim_{n\to\infty} \mu_n  = \infty$, and all 
these eigenvalues have finite multiplicity;
\item[(f)] There exists a complete orthonormal base of $L^2(\bb T^d)$ composed of eigenvectors of $- \mc L_\Lambda$.
\end{itemize}
In view of (a), (c) and (d), by the Hille-Yoshida theorem, $\mc L_\Lambda$ is the generator of a strongly continuous contraction semigroup in $L^2(\bb T^d)$. 
\end{theorem}

The notation $\mc H^1_\Lambda$ will be defined in Section \ref{3}, and it has been chosen 
in analogy to the notation used for Sobolev spaces.

\medskip

\subsection{The hydrodynamic equation}
\label{ss2.3}

Consider a bounded Borel measurable profile $\rho_0 : \bb T^d
\to \bb R$. A bounded function $\rho : \bb R_+\times \bb T^d \to \bb R$
is said to be a weak solution of the parabolic differential equation
\begin{equation}
\label{edp}
\left\{
\begin{array}{l}
{\displaystyle \partial_t \rho \; =\; \mc L_\Lambda \rho } \\
{\displaystyle \rho(0,\cdot) \;=\; \rho_0(\cdot)} \, ,
\end{array}
\right.
\end{equation}
if for all functions  $H$ in $\mc H^1_\Lambda $ and all $t>0$, $\rho$ satisfies the integral equation
\begin{equation}\label{eqint1}
\dl \rho_t, H\dr \;-\; \dl \rho_0 , H\dr 
- \int_0^t \dl \rho_s , \mc L_\Lambda  H \dr \, ds\; \;=\; 0,
\end{equation}
where $ \rho_t$ is the notation for $\rho(t,\cdot)$.
We prove in Subsection \ref{uniqueness} the uniqueness of weak solutions of \ref{eqint}. Existence
follows from the convergence result for the empirical measures associated to the diffusively scaled 
exclusion processes with slow bonds over $\Lambda$, this is discussed in Section \ref{4}.

\begin{theorem}
\label{t01}
Fix a Borel measurable initial profile $\gamma : \bb T^d \to [0,1]$ and
consider a sequence of probability measures $\mu_N$ on $\Omega_N$ associated to $\gamma$. Then, for any $t\geq 0$,
\begin{equation*}
\lim_{N\to\infty}
\bb P^N_{\mu_N} \Big\{ \, \Big\vert \pfrac{1}{N^d} \sum_{x\in\bb T_N^d} 
H(x/N)\, \eta_t(x) - \int_{\bb{T}^d} H(u) \rho(t,u) du \Big\vert 
> \delta \Big\} \;=\; 0 \, ,
\end{equation*}
for every $\delta>0$ and every function $H\in C(\bb{T}^d)$, where $\rho$
is the unique weak solution of the differential equation \eqref{edp} with $\rho_0=\gamma$.
\end{theorem}

\section{The operator $\mc L_\Lambda$}\label{3}
We begin by studying properties of $\mc L_\Lambda$ defined in the domain $\mf D_\Lambda$ and we consider
the extension afterwards.
\begin{lemma}\label{denso}
 The domain $\mf D_\Lambda$ is dense in $L^2(\bb T^d)$. 
\end{lemma}
\begin{proof}
 It is enough to prove that there exists a subset of $\mf D_\Lambda$ which is dense in $L^2(\bb T^d)$. 
All smooth functions with support contained in $\bb T^d\backslash \partial\Lambda$ belong
to $\mf D_\Lambda$, which is clearly a dense subset of $L^2(\bb T^d)$, since $\p\Lambda$ is a smooth zero Lebesgue measure surface that divides $\bb T^d\backslash \partial\Lambda$ in two disjoint open regions.
\end{proof}

\medskip

From now on, we use $\ell_d$ to denote the $d$-dimensional Lebesgue measure on $\bb T^d$.

\medskip

\begin{lemma}\label{simetrico}
 The operator $-\mc L_\Lambda: \mf D_\Lambda\rightarrow L^{2}(\bb T^d)$ is symmetric and non-negative.
 Futhermore, it satisfies a Poincar\'e inequality, which means that there exists a finite constant
$C>0$ such that
\begin{equation}
\Vert H \Vert^2\leq C \,\dl -\mc L_\Lambda H, H\dr + \Big(\int_{\bb T^d} H(x)\,dx\Big)^2
\end{equation}
for all functions $H\in \mf D_\Lambda$.
\end{lemma}
\begin{proof}
Let $H, G\in \mf D_\Lambda$. Write $H=h+\lambda_h  \,{\bf 1}_\Lambda$ and 
$G=g+\lambda_g\,{\bf 1}_\Lambda$, as in definition \ref{operator}. By the
first Green identity and condition (iii) in definition \ref{operator}, we have that
\begin{eqnarray}
\label{green}
\lambda_h \int_{\Lambda}\Delta g\, du & = & \lambda_h\int_{\partial\Lambda} ( \nabla g \cdot \vec{\zeta} ) \ dS
\; = \; - \lambda_h \, \lambda_g \, \textrm{Vol}_{d-1} (\partial\Lambda) \\
& = & \lambda_g \int_{\partial\Lambda} ( \nabla h \cdot \vec{\zeta} ) \ dS
\; = \; \lambda_g \int_{\Lambda}\Delta h\, du \, , \nonumber
\end{eqnarray}
where $dS$ is a infinitesimal element of volume of $\partial\Lambda$ and $\textrm{Vol}_{d-1} (\partial\Lambda)$ is its $(d-1)$-dimensional volume.
Thus,
\begin{eqnarray*}
\dl H,-\mc L_\Lambda G\dr & = & \dl h+\lambda_h\, {\bf 1}_\Lambda,-\Delta g\dr=-\int_{\bb T^d} h\,\Delta g\, du-
\lambda_h\int_{\Lambda}\Delta g\, du \\
& = & - \int_{\bb T^d} g \, \Delta h \, du -
 \lambda_g \int_{\Lambda}\Delta h\, du
\; = \;  \dl -\mc L_\Lambda H,G \dr \,.
\end{eqnarray*}
For the non-negativeness, using \ref{green} above,
\begin{eqnarray*}
 \dl H,-\mc L_\Lambda H\dr & = & -\int_{\bb T^d} h\,\Delta h\, du- \lambda_h\int_{\Lambda}\Delta h \, du \\
 & = & \int_{\bb T^d}\vert \nabla h\vert^2 \, du + \lambda_h^2 \, \textrm{Vol}_{d-1} (\partial\Lambda)
\, \geq \, 0\,.
\end{eqnarray*}
It remains to prove the Poincar\'e inequality. Write
\begin{equation*}
\Vert H \Vert^2 - \Big(\int_{\bb T^d} H(x)\,dx\Big)^2 = \int_{\bb T^d}\Big[H(u)-\int_{\bb T^d} H(v)\; dv\Big]^2 \, du \, ,
\end{equation*}
which can be rewritten as
\begin{equation*}
\int_{\bb T^d} \Big[\Big(h(u)-\int_{\bb T^d}h(v) \, dv\Big) +\lambda_h \Big(\textbf{1}_\Lambda(u)-\ell_d(\Lambda)\Big)
\Big]^2 \, du \, .
\end{equation*}
Now apply the inequality $(a+b)^2\leq 2\,(a^2+b^2)$ to the previous expression to obtain that it is bounded by
\begin{equation*}
2\,\int_{\bb T^d} \Big(h(u)-\int_{\bb T^d}h(v) \, dv\Big)^2 \, du +2\,\lambda_h^2 \Big( \ell_d(\Lambda)-(\ell_d(\Lambda))^2\Big)\,.
\end{equation*}
By the usual Poincar\'e inequality, see \cite{e}, the last expression is less than or equal to
\begin{equation*}
2\,C_1\,\int_{\bb T^d} \vert \nabla h(u)\vert^2 \, du  
+2\,\lambda_h^2 \Big( \ell_d(\Lambda)-(\ell_d(\Lambda))^2\Big)\,.
\end{equation*}
Choosing a constant $C_2>0$ such that $\ell_d(\Lambda)-(\ell_d(\Lambda))^2\leq
\,C_2 \textrm{Vol}_{d-1} (\partial\Lambda)$, the previous expression is bounded above by
$$
2\,\max\{C_1,C_2\}\,\dl-\mc L_\Lambda H, H\dr \, ,
$$
which finishes the proof with $C = 2 \,\max\{C_1,C_2\}$.
\end{proof}

\medskip

Denote by $\dl \cdot, \cdot\dr_{1,\Lambda}$ the inner product on $\mf D_\Lambda$
defined by
\begin{equation*}
\dl F,G\dr_{1,\Lambda} \; =\; \dl F,G\dr \; +\; \dl F,-\,\mc L_\Lambda G\dr\,.
\end{equation*}
Let $\mc H^1_\Lambda $ be the set of all functions $F$ in $L^2(\bb T^d)$
for which there exists a sequence $\{F_n : n\ge 1\}$ in $\mf D_\Lambda$
such that $F_n$ converges to $F$ in $L^2(\bb T^d)$ and $F_n$ is Cauchy
for the inner product $\dl \cdot, \cdot \dr_{1,\Lambda}$. Such sequence
$\{F_n\}$ is called admissible for $F$. For $F$, $G$ in
$\mc H^1_\Lambda $, define
\begin{equation}
\label{f20}
\dl F,G\dr_{1,\Lambda}\; =\; \lim_{n\to\infty} \dl F_n,G_n\dr_{1,\Lambda}\,,
\end{equation}
where $\{F_n\}$, $\{G_n\}$ are admissible sequences for $F$, $G$,
respectively.  By \cite[Proposition 5.3.3]{z}, the limit exists and
does not depend on the admissible sequence chosen. Moreover, $\mc H^1_\Lambda $
 endowed with the scalar product $\dl \cdot, \cdot \dr_{1,\Lambda}$
just defined is a real Hilbert space. From now on, we consider $\mc H^1_\Lambda $
with the norm induced by $\dl \cdot, \cdot \dr_{1,\Lambda}$, unless we 
mention that we are going to use the $L^2$-norm.

\begin{lemma}\label{embedding}
The embedding $\mc H^1_\Lambda  \subset L^2(\bb T^d)$ is compact.
\end{lemma}
\begin{proof}
Let $\{ H_n\}$ a bounded sequence in $\mc H^1_\Lambda $. Fix $\{ F_n\}$ as a sequence in $\mf D_\Lambda$
such that $\Vert F_n-H_n\Vert \to 0$ and $\{ F_n\}$ is also bounded in
$\mc H^1_\Lambda $. Thus, to get a convergent subsequence of $\{ H_n\}$,
 it is sufficient to find a convergent subsequence of $\{F_n\}$ in  $L^2(\bb T^d)$. Write
$F_n=f_n+\lambda_n\textbf{1}_\Lambda$, with $f_n\in C^2(\bb T^d)$. 
Then,
\begin{equation*}
 \dl F_n,F_n\dr_{1,\Lambda}=\dl f_n+\lambda_n\textbf{1}_\Lambda,f_n+\lambda_n\textbf{1}_\Lambda\dr+
\dl f_n+\lambda_n\textbf{1}_\Lambda,-\Delta f_n \dr \, .
\end{equation*}
Expanding the right hand side and using \ref{green}, we get that
\begin{equation*}
\dl F_n,F_n\dr_{1,\Lambda}= \Vert f_n\Vert^2+\lambda^2_n\ell_d(\Lambda)+2\lambda_n\int_{\Lambda}f_n(u)\,du
+ \Vert \nabla f_n \Vert^2+ \lambda_n^2 \, \textrm{Vol}_{d-1} (\partial\Lambda)\,,
\end{equation*}
which is greater or equal to
\begin{equation*}
 \Vert f_n\Vert^2+\lambda^2_n\ell_d(\Lambda)- \lambda_n^2-
\ell_d(\Lambda)\int_{\Lambda}f_n^2(u)\,du
+ \Vert \nabla f_n \Vert^2+\lambda_n^2 \, \textrm{Vol}_{d-1} (\partial\Lambda)
\end{equation*}
\begin{equation*}
=\Big(\ell_d(\Lambda)-1+ \textrm{Vol}_{d-1} (\partial\Lambda) \Big)\,
\lambda^2_n + (1-\ell_d(\Lambda))\,\int_{\Lambda}f_n^2(u)\,du +
\int_{\Lambda^\complement}f_n^2(u)\,du + \Vert \nabla f_n \Vert^2
\end{equation*}
\begin{equation*}
\geq \Big( \textrm{Vol}_{d-1} (\partial\Lambda) - \ell_d(\Lambda^\complement) \Big)\,
\lambda^2_n + (1-\ell_d(\Lambda))\,\Vert f_n \Vert^2
 + \Vert \nabla f_n \Vert^2\,.
\end{equation*}
If we put $\tilde{f}_n = f_n+1$, and write $F_n = \tilde{f}_n - \lambda_n \textbf{1}_{\Lambda^\complement}$, an analogous computation shows that $\dl F_n,F_n\dr_{1,\Lambda}$ is greater or equal than
\begin{equation*}
 \Big( \textrm{Vol}_{d-1} (\partial\Lambda) - \ell_d(\Lambda) \Big)\,
\lambda^2_n + (1-\ell_d(\Lambda^\complement))\,\Vert f_n \Vert^2
 + \Vert \nabla f_n \Vert^2\,.
\end{equation*}
By the classical isoperimetric inequality on the Torus (see \cite[Lemma 4.6]{ct} for the statement and a direct proof), we have that
$$
\max\{ \, \textrm{Vol}_{d-1} (\partial\Lambda) - \ell_d(\Lambda^\complement) \, , \, \textrm{Vol}_{d-1} (\partial\Lambda) - \ell_d(\Lambda) \, \} > 0 \, .
$$
Since $\{\dl F_n,F_n\dr_{1,\Lambda}\}$ is a bounded sequence, we conclude that 
$\{\lambda_n\}$ is bounded, as well
the sequence $\{ \Vert f_n\Vert^2 + \Vert \nabla f_n\Vert^2\}$.
By the Rellich-Kondrachov Compactness Theorem, see \cite[Theorem 5.7.1]{e}, $\{f_n\}$ 
has a convergent subsequence in $L^2(\bb T^d)$. From this subsequence, choosing a 
convergent subsequence of $\{\lambda_n\}$ finishes the proof. 
\end{proof}

\medskip

\begin{lemma}\label{image}
The image of $\bb I-\, \mc L_\Lambda:\mf D_\Lambda\to L^2(\bb T^d)$ is dense in
$L^2(\bb T^d)$.
\end{lemma}
\begin{proof}
By a similar argument to the one found in Lemma \ref{denso}, it is enough to show that any smooth function $f$
with support contained in $\bb T^d\backslash \partial\Lambda$ belongs to 
$(\bb I- \, \mc L_\Lambda)(\mf D_\Lambda)$. Therefore, we need to find a function $h$ in $C^2(\bb T^d)$ with support in $\bb T^d\backslash\partial\Lambda$ such that 
\begin{equation*}
h- \Delta h \,= \,f\,.
\end{equation*}
From the classical theory of second-order parabolic equations,  e.g., see \cite[Theorem 5.7.1]{e}, this equation has a smooth solution, which finishes the proof.

\end{proof}

\begin{proof}\textbf{(Proof of Theorem \ref{Lambda})}
\smallskip

\textbf{(a)} Since $\mf D_\Lambda \subset \mc H^1_\Lambda $, it follows from Lemma \ref{denso} that $\mc H^1_\Lambda $ is dense in $L^2(\bb T^d)$. 
\smallskip

\textbf{(b)} Denote $\bb I-\,\mc L_\Lambda = \mc A : \mf D_\Lambda \to
\bb L^2(\bb T^d)$. From Lemma \ref{simetrico}, $\mc A$ is linear, symmetric and strongly monotone on the
Hilbert space $L^2(\bb T^d)$. By strongly monotone, we mean that there exists $c>0$ such that
\begin{equation*}
 \dl \mc A \,H, H\dr \geq \, c \, \Vert H\Vert^2\,,\quad \forall H\in \mf D_\Lambda\,.
\end{equation*}
In this case, $\mc A$ satisfy the inequality above with $c=1$.
By \cite[Theorem 5.5.a]{z}, in the conditions above, the Friedrichs extension 
$\mc A:\mc H^1_\Lambda \to L^2(\bb T^2)$ is self-adjoint, bijective and strongly monotone. 
By an abuse of notation, define now the extension $\mc L_\Lambda:\mc H^1_\Lambda \to L^2(\bb T^2)$ 
as $(\bb I-\mc A)$.
Since $\bb I$ and $\mc A$ are self-adjoint in $\mc H^1_\Lambda $, this property is 
inherited by $\mc L_\Lambda:\mc H^1_\Lambda \to L^2(\bb T^2)$. 

For non-positiveness, note that
\begin{equation*}
\dl -\,\mc L_\Lambda\, H,H\dr=\dl-(\bb I-\mc A)H,H\dr=-\dl H,H\dr + \dl \mc A\, H,H\dr \geq 0\,.
\end{equation*}

\textbf{(c)} As mentioned in the proof of (b) above, the Friedrichs extension 
$\mc A:\mc H^1_\Lambda \to L^2(\bb T^2)$ is bijective. So it remains to show that $\mf D_\Lambda$ is a core of $\mc A:\mc H^1_\Lambda \to L^2(\bb T^2)$. For any operator $B$, denote by $\mc G(B)$ the graphic of $B$. Then $\mf D_\Lambda$ is a core for $\mc A$, if the closure of $\mc G(\mc A|_{\mf D_\Lambda})^{_{L^2\times L^2}}$ in $L^2\times L^2$ is equal to $\mc G(\mc A)$. Since $\mc A$ is self-adjoint, $\mc A$ is a closed operator, or else, $\mc G(\mc A)$ is a closed set. Thus the closure of
$\mc G(\mc A|_{\mf D_\Lambda})$ is a subset of $\mc G(\mc A)$.
Let $H\in\mc H^1_\Lambda$, from Lemma \ref{image}, there exists
a sequence $\{H_n\}$ in $\mf D_\Lambda$ such that $\mc A\, H_n$ converges to $\mc A\, H$ in $L^2$.
Hence, as proved in  \cite[Theorem 5.5.a]{z}, $\mc A^{-1}$ is a bounded linear operator, 
and $H_n$ converges to $H$ in $L ^2$, which yields that the closure of
$\mc G(\mc A|_{\mf D_\Lambda})$ contains $\mc G(\mc A)$.
\smallskip

\textbf{(d)} Fix a function $H$ in $\mc H^1_\Lambda $ and $\mu >0$. Put $G
= (\mu \bb I - \mc L_\Lambda) H$. Taking the inner product with respect
to $H$ on both sides of this equality, we obtain that
\begin{eqnarray*}
\mu\, \dl H , H\dr \;+\; \dl -\mc L_\Lambda H , H\dr
\;=\; \dl H , G\dr \;\le\;  \dl H , H\dr^{1/2} \,  
\dl G , G\dr^{1/2}\;.
\end{eqnarray*}
Since $H$ belongs to $\mc H^1_\Lambda $, by (b), the second term on the left
hand side is positive. Therefore, $\mu\Vert  H\Vert \le \Vert G \Vert =
\Vert (\mu \bb I - \mc L_\Lambda) H \Vert $.
\smallskip

\textbf{(e) and (f)} We have seen that the operator $(\bb I - \,\mc L_\Lambda): \mf D_\Lambda \to
L^2(\bb T)$ is symmetric and strongly monotone.  By Lemma \ref{embedding} ,
the embedding $\mc H^1_\Lambda  \subset L^2(\bb T^d)$ is compact. Therefore,
by \cite[Theorem 5.5.c]{z}, the Friedrichs extension $\mc A : \mc H^1_\Lambda  \to L^2(\bb T^d)$, satisfies claims
(e) and (f) with $1\le \lambda_1 \le \lambda_2 \le \cdots$,
$\lambda_n\uparrow\infty$. In particular, the operator $- \mc L_\Lambda =
(\mc A - \bb I)$ has the same property with $0\le \mu_1 \le
 \mu_2 \le \cdots$, $ \mu_n\uparrow\infty$. Since $0$ is an
eigenvalue of $-\mc L_\Lambda$, a constant function is an eigenfunction with eigenvalue $0$, then (e) and
(f) also hold.
\end{proof}

\section{Scaling Limit}\label{4}

Let $\mc M$ be the space of positive Radon measures on $\bb T^d$ with total
mass bounded by one endowed with the weak topology. For a measure $\pi \in \mc M$ 
and a measurable $\pi$-integrable function $H:\bb T^d \to \bb R$, 
we denote by $\<\pi , H\>$ the integral of $H$ with respect to $\pi$.

Recall that $\{\eta^N_t : t\ge 0\}$ denote a Markov process with state space $\Omega_N$ and generator 
$L_N$ speeded up by $N^2$. Let $\pi^{N}_{t} \in \mc M$ be the empirical measure at time $t$ associated 
to $\{\eta^N_t : t\ge 0\}$, which is the random measure in $\mc M$ given by 
\begin{equation}
\label{f01}
\pi^{N}_{t} \;=\; \frac{1}{N^d} \sum _{x\in \bb T_N^d} \eta^N_t (x)\,
\delta_{x/N}\;,
\end{equation}
where $\delta_u$ is the Dirac measure concentrated on $u$. 

Note that
\begin{equation*}
\<\pi^N_t, H\> \;=\; \pfrac 1{N^d} \sum_{x\in\bb T_N^d}
H (\pfrac xN) \eta^N_t(x)\; ,
\end{equation*}
for the empirical measures, and $\<\pi , H\> = \dl \rho , H \dr$, for absolutely continuous
measures $\pi$ with $L^2$ bounded density $\rho$, and $H \in L^2 (\bb T^d)$.

Fix $T>0$. Let $D([0,T], \mc M)$ be the space of $\mc M$-valued
c\`adl\`ag trajectories $\pi:[0,T]\to\mc M$ endowed with the
\emph{Skorohod} topology. Then, the $\mc M$-valued process $\{\pi^N_t:t\ge 0\}$
is a random element of $D([0,T], \mc M)$ whose distribution is determined by
the initial distribution of $\{\eta^N_t : t\ge 0\}$. For each probability measure $\mu$ on
$\Omega_N$, denote by $\bb Q_{\mu}^{\Lambda,N}$ the distribution of $\{\pi^N_t:t\ge 0\}$
on the path space $D([0,T], \mc M)$, when $\eta^N_0$ has distribution $\mu$. 

For a Borel measurable profile $\gamma : \bb T^d \to [0,1]$, suppose that there exists
a unique weak solution $\rho$ of \eqref{edp} with initial condition $\gamma$.
We denote by $\bb Q_\Lambda^\gamma$ be the probability measure on $D([0,T], \mc M)$ concentrated 
on the deterministic path  $\pi(t,du) = \rho (t,u)du$.

\begin{proposition}
\label{s15}
Fix a Borel measurable profile $\gamma : \bb T^d \to [0,1]$ and consider a
sequence $\{\mu_N : N\ge 1\}$ of measures on $\Omega_N$
associated to $\gamma$ in the sense of (\ref{f09}). Then there exists
a unique weak solution $\rho$ of \eqref{edp} with initial condition $\gamma$
and the sequence of probability measures $\bb Q_{\mu_N}^{\Lambda,N}$ 
converges weakly to $\bb Q_\Lambda^\gamma$ as $N\uparrow\infty$.
\end{proposition}

It is straightforward to obtain Theorem \ref{t01} as a corollary of the previous proposition.
The proof of Proposition \ref{s15} follows directly from the uniqueness of weak solutions
of \eqref{edp} and the next two results:

\smallskip

\begin{proposition}
\label{s06}
For any sequence $\{\mu_N : N\ge 1\}$ of probability measures with $\mu_N$ concentrated on $\Omega_N$, 
the sequence of measures $\{\bb Q_{\mu_N}^{\Lambda,N} : N\ge 1\}$ is tight.
\end{proposition}

\begin{proposition}\label{3.5} 
Fix a Borel measurable profile $\gamma : \bb T^d \to [0,1]$ and consider a sequence $\{\mu_N : N\ge 1\}$ of probability 
measures on $\Omega_N$ associated to $\gamma$ in the sense of (\ref{f09}). Then any limit point of 
$\bb Q_{\mu_N}^{\Lambda,N}$ is concentrated in absolutely continuous trajectories that are weak solutions of \eqref{edp} with initial condition $\gamma$.
\end{proposition}
\smallskip

\begin{proof} {\emph{(of Proposition \ref{s15}).}}
By Proposition \ref{s06}, $\bb Q_{\mu_N}^{\Lambda,N}$ is tight. Thus, if a limit point $\bb Q^*$ is concentrated in 
solutions of equation \eqref{edp}, by the uniqueness result proved  in section \ref{uniqueness}, we have that 
$\bb Q^* = \bb Q^\gamma_\Lambda$. And the statement follows from Proposition \ref{3.5}.
\end{proof}

\smallskip

In Subsection \ref{tightness}, we prove Proposition \ref{s06} and in Subsection \ref{limitpoints} we 
show Proposition \ref{3.5}. As a consequence, we have the existence of solutions of \eqref{edp} with initial 
condition $\gamma$. We complete the proof in subsection \eqref{uniqueness} showing the uniqueness of weak 
solutions of \eqref{edp}.

\subsection{Tightness}\label{tightness}

Here we prove Proposition \ref{s06}.
Let $D([0,T], \bb R)$ be the space of $\bb R$-valued c\`adl\`ag trajectories with domain $[0,T]$ endowed with the
\emph{Skorohod} topology. To prove tightness of $\{\pi^{N}_t : 0\le t \le T\}$ 
in $D([0,T], \mc M)$, it is enough to show tightness in $D([0,T], \bb R)$ of the real-valued processes 
$\{\<\pi^{N}_t ,H\> : 0\le t \le T\}$ for a set of functions $H:\bb T^d\to \bb R$ which is dense in the space 
of continuous real functions on $\bb T^d$ endowed with the uniform topology, see \cite{kl}. Futhermore, if a sequence of distributions in $D([0,T], \bb R)$ endowed with the uniform topology is tight, then it is also tight in $D([0,T], \bb R)$ 
endowed with the Skorohod topology. Here we prove tightness of $\{\<\pi^{N}_t ,H\> : 0\le t \le T\}$ in $D([0,T], \bb R)$, endowed with the uniform topology, for $H\in C^2(\bb T^d)$. 
\smallskip

Fix $H\in C^2(\bb T^d)$. By definition $\{\<\pi^{N}_t ,H\> : 0\le t \le T\}$ is tight in $D([0,T], \bb R)$ endowed with the uniform topology if, for the boundedness, 
\begin{equation}
\label{tight1}
\lim_{m \rightarrow \infty} \, \sup_{N} \, \bb P^N_{\mu_N}\left[ \sup_{0\le t \le T} | \<\pi^{N}_t ,H\> | > m \right] =  0\,,
\end{equation}
and, for the equicontinuity, 
\begin{equation}
\label{tight2}
\lim_{\delta \rightarrow 0} \, \limsup_{N \rightarrow \infty} \, \bb P^N_{\mu_N}\left[ \sup_{|t-s| \le \delta} | \<\pi^{N}_t ,H\> -  \<\pi^{N}_s ,H\>| > \epsilon \right] =  0\,, \, \textrm{ for all } \epsilon >0 \, .
\end{equation}
The limit in (\ref{tight1}) is trivial since
$$
|\<\pi^{N}_t ,H\>| \le  \sup_{0\le t \le T} |H(t)| \, .
$$
So we only need to prove (\ref{tight2}). By Dynkyn's formula (see appendix in \cite{kl}),
\begin{equation}\label{M}
M^{N}_{t}=\<\pi^{N}_{t}, H\>- \<\pi^{N}_{0}, H\>-\int_{0}^{t}N^{2}L_{N}\<\pi^{N}_{s},H\>ds
\end{equation}
is a martingale. By the previous expression, (\ref{tight2}) follows from
\begin{equation}
\label{tight3}
\lim_{\delta \rightarrow 0} \, \limsup_{N \rightarrow \infty} \, \bb P^N_{\mu_N}\left[ \sup_{|t-s| \le \delta} | M^{N}_{t} -  M^{N}_{s} | > \epsilon \right] =  0\,, \, \textrm{ for all } \epsilon >0 \, ,
\end{equation}
and
\begin{equation}
\label{tight4}
\lim_{\delta \rightarrow 0} \, \limsup_{N \rightarrow \infty} \, \bb P^N_{\mu_N}\left[ \sup_{0 \le t-s \le \delta} \Big| \int_{s}^{t}N^{2}L_{N}\<\pi^{N}_{s},H\>ds \Big| > \epsilon \right] =  0\,, \, \textrm{ for all } \epsilon >0 \, .
\end{equation}
Indeed, we show the stronger results below:
\begin{equation}
\label{tight5}
\lim_{\delta \rightarrow 0} \, \limsup_{N \rightarrow \infty} \, \bb E^N_{\mu^N} \left[ \sup_{|t-s| \le \delta} | M^{N}_{t} -  M^{N}_{s} | \right] =  0 \, ,
\end{equation}
and
\begin{equation}
\label{tight6}
\lim_{\delta \rightarrow 0} \, \limsup_{N \rightarrow \infty} \, \bb E^N_{\mu^N} \left[ \sup_{0 \le t-s \le \delta} \Big| \int_{s}^{t}N^{2}L_{N}\<\pi^{N}_{s},H\>ds \Big| \right] =  0 \, .
\end{equation}
To verify (\ref{tight5}), we use the quadratic variation of $M^{N}_{t}$ that we denote by $\<M^{N}_{t}\>$.
By Doob's inequality, we have that
\begin{eqnarray*}
\bb E^N_{\mu^N} \left[ \sup_{|t-s| \le \delta} | M^{N}_{t} -  M^{N}_{s} | \right] & \le &
2  \, \bb E^N_{\mu^N} \left[ \sup_{0 \le t \le T} |M^{N}_{t}| \right] \\ 
& \le & 2 \, \bb E^N_{\mu^N} \left[ \sup_{0 \le t \le T} |M^{N}_{t}|^2 \right]^{\pfrac{1}{2}} \, \le \,
2 \, \bb E^N_{\mu^N} \left[ \sup_{0 \le t \le T} \<M^{N}_{t}\> \right]^{\pfrac{1}{2}}.
\end{eqnarray*}
Since
$$
\<M^{N}_{t}\>=\int_{0}^{t} N^{2} [ L_{N}\<\pi^{N}_{s},H\>^{2}-2\<\pi^{N}_{s},H\>L_{N}\<\pi^{N}_{s},H\>] ds \, ,
$$
we obtain by a straightforward computation that 
\begin{equation*}
\<M^{N}_{t}\> = \int_{0}^{t} N^{2} \sum_{j=1}^d\sum_{x\in\bb T_{N}^d} \xi^{N}_{x,x+e_j}\pfrac{1}{N^{2d}}\Big[(\eta_{s}(x)-\eta_{s}(x+e_j))
(H(\pfrac{x+e_j}{N})-H(\pfrac{x}{N}))\Big]^{2} ds \, .
\end{equation*}
Therefore, since $\xi^N_{x,x+e_j}\leq 1$, 
\begin{eqnarray}\label{ptos lim}
\<M^{N}_{t}\> & \leq & \frac{T}{N^{2d-2}}\sum_{j=1}^d\sum_{x\in\bb T_{N}^d} \xi^{N}_{x,x+e_j} 
\Big[H(\pfrac{x+e_j}{N})-H(\pfrac{x}{N})\Big]^2 \nonumber \\
& \leq & \frac{Td}{N^{d}} \, \Big( \sup_{0 \le t \le T} |\nabla H (t) \cdot e_j| \Big)^2\,. 
\end{eqnarray}
Thus, $M^{N}_{t}$ converges to zero in $L^{2}$ and (\ref{tight5}) holds.

We finish the proof by verifying (\ref{tight6}). Write
\begin{equation*}
N^{2}L_{N}\<\pi^{N}_{s},H\> = \pfrac{1}{N^{d-2}}\sum_{j=1}^d\sum_{x\in \bb T_{N}^d}\xi^{N}_{x,x+e_j}((\eta_{s}(x)-\eta_{s}(x+e_j))\left(H(\pfrac{x+e_j}{N})
-H(\pfrac{x}{N})\right)
\end{equation*}
\begin{equation*}
= \pfrac{1}{N^{d-2}}\sum_{j=1}^d\sum_{x\in \bb T_{N}^d}\eta_{s}(x)
\Big[\xi^{N}_{x,x+e_j}\Big(H(\pfrac{x+e_j}{N})-H(\pfrac{x}{N})\Big)+\xi^{N}_{x,x-e_j}
\Big(H(\pfrac{x-e_j}{N})-H(\pfrac{x}{N})\Big)\Big]\,.
\end{equation*}
Define $\Gamma_N\subset\bb T^d_N$ as the set of vertices
whose have some adjacent edge with exchange rate not equal to one. Then $N^{2}L_{N}\<\pi^{N}_{s},H\>$ is equal to
\begin{equation*}
\pfrac{1}{N^{d-2}}\sum_{j=1}^d\sum_{x\notin\Gamma_N}
\eta_{s}(x)\Big[H(\pfrac{x+e_j}{N})+H(\pfrac{x-e_j}{N})-2H(\pfrac{x}{N})\Big]
\end{equation*}
\begin{equation*}
+\pfrac{1}{N^{d-2}}\sum_{j=1}^d\sum_{x\in\Gamma_N}
\eta_{s}(x)\Big[\xi^{N}_{x,x+e_j}\Big(H(\pfrac{x+e_j}{N})-H(\pfrac{x}{N})\Big)+\xi^{N}_{x,x-e_j}
\Big(H(\pfrac{x-e_j}{N})-H(\pfrac{x}{N})\Big)\Big] \, .
\end{equation*}
By the Taylor expansion (remember  $H\in C^2$), the absolute value of the first term above is bounded by  $\sup_{0\le t \le T} | \Delta H (t)|$. Since there are in order of $N^{d-1}$ terms in $\Gamma_N$, and
$\xi_{x,x+e_j}\leq 1$,  the absolute value of the
second term above is bounded by 
\begin{equation*}
\pfrac{1}{N^{d-2}}\sum_{j=1}^d\sum_{x\in\Gamma_N}
\vert H(\pfrac{x+e_j}{N})-H(\pfrac{x}{N})\vert+\vert H(\pfrac{x-e_j}{N})-H(\pfrac{x}{N})\vert
\leq 2d\, \sup_{0 \le t \le T} |\nabla H (t) \cdot e_j| \,.
\end{equation*}
Therefore, there exists $C>0$, depending only on $H$, such that
$|N^{2}L_{N}\<\pi^{N}_{s},H\>| \leq C$, which yields
\begin{equation*}
 \left|\int_{r}^{t}N^{2}L_{N}\<\pi^{N}_{s},H\>ds\right|\leq C(t-r)\,.
\end{equation*}
and (\ref{tight6}) holds.  

\smallskip

\subsection{Caracterization of limit points} 
\label{limitpoints} 
Let $\gamma : \bb T^d \to [0,1]$ be a Borel measurable profile and consider a
sequence $\{\mu_N : N\ge 1\}$ of measures on $\Omega_N$
associated to $\gamma$ in the sense of (\ref{f09}).
We prove Proposition \ref{3.5} in this subsection, i.e., that all limit points $\bb Q^*$ of the
sequence $\bb Q_{\mu_N}^{\Lambda,N}$ are concentrated on absolutely
continuous trajectories $\pi(t,du) = \rho(t,u) du$, whose density
$\rho(t,u)$ is a weak solution of the hydrodynamic equation
\eqref{edp} with $\gamma$ as the initial condition.

Let $\bb Q^*$ be a limit point of the sequence $\bb Q_{\mu_N}^{\Lambda,N}$
and assume, without loss of generality, that $\bb Q_{\mu_N}^{\Lambda,N}$
converges to $\bb Q^*$.

Since there is at most one particle per site, $\bb
Q^*$ is concentrated on trajectories $\pi_t(du)$ which are absolutely
continuous with respect to the Lebesgue measure, $\pi_t(du) =
\rho(t,u) du$, and whose density $\rho$ is non-negative and bounded by
1, se \cite[Chapter 4]{kl}. 

\smallskip

We shall prove the following result:

\begin{lemma}\label{c1}
Any limit point $\bb Q^*$ of $\bb Q_{\mu_N}^{\Lambda,N}$ is concentrated is absolutely continuous trajectories
$\pi_t(du) = \rho(t,u) du$ such that, for any $H\in\mf D_\Lambda$, 
\begin{equation}\label{equation}
\dl \rho_t,  H \dr - \dl \gamma , H \dr \;=\; \int_0^t  \, \dl \rho_s \,,\, \mc L_\Lambda  H \dr \;ds\,.
\end{equation}
\end{lemma}

\smallskip

With the previous lemma we can show Proposition \ref{3.5}.

\begin{proof} {\emph{(of Proposition \ref{3.5}).}} 
It just remains to extend the equality \eqref{equation} to functions $H\in \mc H^1_\Lambda $.
Let $H\in\mc H^1_\Lambda $. Since $\bb I -\kappa\,\mc L_\Lambda: \mf D_\Lambda\to L^2(\bb T^d)$ 
is a core for the Friedrichs extension, there exists a sequence $H_n\in\mf D_\Lambda$ such that
\begin{equation*}
(H_n, (\bb I -\kappa\,\mc L_\Lambda)\,H_n)\to (H, (\bb I -\kappa\,\mc L_\Lambda)\,H)
\end{equation*}
in $L^2(\bb T^d)\times L^2(\bb T^d)$. Thus, $H_n\to H$ and $\mc L_\Lambda\,H_n\to\mc L_\Lambda\,H$, both
in $L^2(\bb T^d)$. Replacing $H_n$ in equality 
\eqref{equation}, and taking the limit as $n \rightarrow \infty$, it finishes the proof.
\end{proof}

\smallskip

The remain of this section is devoted to the proof of Lemma \ref{c1}. Fix a function $H\in \mf D_\Lambda$ and define the martingale
$M^{N}_t$ by
\begin{equation}\label{M1}
\<\pi^N_t, H\> \,-\, \<\pi^N_0, H \> \,-\,
\int_0^t  \, N^2 L_N \<\pi^N_s ,H \> \,ds\,.
\end{equation}
We claim that, for every $\delta>0$,
\begin{equation}\label{limprob1}
\lim_{N\to\infty} \bb P^N_{\mu_N} \Big[ \sup_{0\le t\le T}
\Big\vert M^{N}_t \Big\vert \, > \, \delta \Big] 
\;=\; 0\,.
\end{equation}
For $H\in C^{2}$, this follows from Chebyshev inequality and the estimates done in the proof of tightness, where we have shown that
\begin{equation}
\label{doob}
\lim_{N\to\infty} \bb E^N_{\mu} \left[ \sup_{0 \le t \le T} |M^{N}_{t}| \right] \, \le \, \lim_{N\to\infty} \bb E^N_{\mu} \left[ \sup_{0 \le t \le T} \<M^{N}_{t}\> \right]^{\pfrac{1}{2}} = 0\, .
\end{equation}
For $H = h+\lambda  \,{\bf 1}_\Lambda$ in $\mf D_\Lambda$, the first inequality in \eqref{ptos lim} is still valid and  
\begin{eqnarray}
\<M^{N}_{t}\> & \leq &
\pfrac{T}{N^{2d-2}}\sum_{j=1}^{d}\sum_{x\in\bb T_{N}^d} \xi^{N}_{x,x+e_j}
\Big[H(\pfrac{x+e_j}{N})-H(\pfrac{x}{N})\Big]^{2} \nonumber \\
\label{sum1}
& = &\pfrac{T}{N^{2d-2}}\sum_{j=1}^{d}\sum_{x\notin\Gamma_N} 
\Big[h(\pfrac{x+e_j}{N})-h(\pfrac{x}{N})\Big]^{2} \\
\label{sum2}
&  &  + \pfrac{T}{N^{2d-2}}\sum_{j=1}^{d}\sum_{x\in\Gamma_N} \xi^{N}_{x,x+e_j}
\Big[H(\pfrac{x+e_j}{N})-H(\pfrac{x}{N})\Big]^{2} \, ,
\end{eqnarray}
where $\Gamma_N$ is also defined in the proof of tightness.
The expression \eqref{sum1} goes to zero as $N$ increases, since the function $h$ is Lipschitz.
For the expression in \eqref{sum2}, let $x\in\Gamma_N$. If $\pfrac{x}{N}\in\Lambda$ and $\pfrac{x+e_j}{N}
\in \Lambda^\complement$, then $\xi^N_{x,x+e_j}\leq \pfrac{1}{N}$. The same occurs if
$\pfrac{x}{N}\in\Lambda^\complement$ and $\pfrac{x+e_j}{N}\in \Lambda$. If $\pfrac{x}{N},\pfrac{x+e_j}{N}
$ both belong to $\Lambda$ or $\Lambda^\complement$, the exchange rate $\xi^N_{x,x+e_j}$ is one, but
$\vert H(\pfrac{x+e_j}{N})-H(\pfrac{x}{N})\vert=\vert h(\pfrac{x+e_j}{N})-h(\pfrac{x}{N})\vert\leq
\pfrac{1}{N} \sup_{0 \le t \le T} |\nabla H (t) \cdot e_j| $. In both cases, the expression \eqref{sum2} is of order $O(N^{-d})$. Therefore, from (\ref{doob}), we obtain \eqref{limprob1}. 

The next step is to show that we can replace $N^2 \bb L_N$ by the continuous operator $\mc L_\Lambda$ in the martingale 
formula (\ref{M1}) and that the resulting expression still converges to zero in probability. This will follow from the ensuing proposition:

\smallskip

\begin{proposition}\label{condicoesI} For any $H\in \mf D_\Lambda$, 
\begin{equation}\label{9}
\lim_{N \rightarrow \infty}
\frac{1}{N^d}\sum_{x\in \bb T_{N}^d}\Big\vert N^{2}\bb{L}_{N}H(\pfrac{x}{N})-\mc L_\Lambda H(\pfrac{x}{N})
\Big\vert = 0\,. 
\end{equation}
\end{proposition}

\begin{proof} As usual, put $H = h+\lambda  \,{\bf 1}_\Lambda$, where $h \in C^2(\bb T^d)$.
Rewrite the sum in \eqref{9} as
\begin{equation*}
\frac{1}{N^d} \sum_{x\notin \Gamma_N}\Big\vert N^{2}\bb{L}_{N}H(\pfrac{x}{N})-\mc L_\Lambda H(\pfrac{x}{N})
\Big\vert+ \frac{1}{N^d}\sum_{x\in \Gamma_N}\Big\vert N^{2}\bb{L}_{N}H(\pfrac{x}{N})-\mc L_\Lambda H(\pfrac{x}{N})\Big\vert\,.
\end{equation*}
The first term above is equal to
\begin{equation*}
\frac{1}{N^d}\sum_{x\notin \Gamma_N}
\Big\vert N^{2}\Big(h(\pfrac{x+e_j}{N})+h(\pfrac{x-e_j}{N})-2h(\pfrac{x}{N})\Big)
-\Delta h(\pfrac{x}{N})\Big\vert\,,
\end{equation*}
which converges to zero because $h\in C^2$. The second one is less than or equal to
the sum of 
\begin{equation}\label{Delta}
\frac{1}{N^d}\sum_{x\in \Gamma_N}\vert\Delta h(\pfrac{x}{N})\vert
\end{equation}
and
\begin{eqnarray}
\label{inGamma} 
\lefteqn{ \!\!\!\!\!\!\!\!\!\!\!\!\!\!\!\!\!\!\!\!\!\!\! \frac{1}{N^{d-1}} \sum_{x\in \Gamma_N}
\sum_{j=1}^d \Big\vert N\xi^N_{x,x+e_j}(H(\pfrac{x+e_j}{N})-H(\pfrac{x}{N})) } \nonumber \\
& & \qquad \qquad + N\xi^N_{x,x-e_j}(H(\pfrac{x-e_j}{N})-H(\pfrac{x}{N})) \Big\vert \, .
\end{eqnarray}
Since there are $O(N^{d-1})$ terms in $\Gamma_N$, the expression in \eqref{Delta}
converges to zero as $N \rightarrow \infty$. Since $\p\Lambda$ is smooth, the quantity of points $x \in \Gamma_N$ 
for which  both $\xi^N_{x,x+e_j}$ and $\xi^N_{x,x-e_j}$ are different of one is negligible. Therefore, we must 
only worry about points $x\in\Gamma_N$ such that, for some $j$, only one of $\xi^N_{x,x+e_j}$ and $\xi^N_{x,x-e_j}$ 
is equal to $N^{-1}$. This occurs in one of the following four cases:  $\pfrac{x}{N}\in \Lambda$, $\pfrac{x-e_j}{N}\in \Lambda$ and $\pfrac{x+e_j}{N}\in \Lambda^\complement$;  $\pfrac{x}{N}\in \Lambda$, $\pfrac{x-e_j}{N}\in \Lambda^\complement$ and $\pfrac{x+e_j}{N}\in \Lambda$;  $\pfrac{x}{N}\in \Lambda^\complement$, $\pfrac{x-e_j}{N}\in \Lambda$ and $\pfrac{x+e_j}{N}\in \Lambda^\complement$;  $\pfrac{x}{N}\in \Lambda^\complement$, $\pfrac{x-e_j}{N}\in \Lambda^\complement$ and $\pfrac{x+e_j}{N}\in \Lambda$. The analysis of these cases are analogous, thus we only consider the first one. Suppose $\pfrac{x}{N}\in \Lambda$, $\pfrac{x-e_j}{N}\in \Lambda$ and $\pfrac{x+e_j}{N}\in \Lambda^\complement$. In this case, 
the summand in \eqref{inGamma} can be rewritten as 
\begin{equation*}
N\xi^N_{x,x+e_j}(H(\pfrac{x+e_j}{N})-H(\pfrac{x}{N}))
+N\xi^N_{x,x-e_j}(H(\pfrac{x-e_j}{N})-H(\pfrac{x}{N})) 
\end{equation*}
\begin{equation*}
 = \vert\vec{\zeta}_{x,j}\cdot e_j\vert\,[H(\pfrac{x+e_j}{N})-H(\pfrac{x}{N})]\
+N\,[H(\pfrac{x-e_j}{N})-H(\pfrac{x}{N})]\,,
\end{equation*}
which becomes uniformly (in $x\in \Gamma_N$) close to
\begin{equation*}
 -\lambda\vert\vec{\zeta}_{x,j}\cdot e_j\vert\,\,
\textrm{sgn}\Big(\vec{\zeta}_{x,j}\cdot e_j\Big)-\pfrac{\p h}{\p u_j}(\pfrac{x}{N})
 =-\lambda \,\vec{\zeta}_{x,j}\cdot e_j-\pfrac{\p h}{\p u_j}(\pfrac{x}{N})\,.
\end{equation*}
The condition $\nabla h|_{\p\Lambda}(u)=-\lambda\, \vec{\zeta}(u)$, which was imposed in the definition
of $\mf D_\Lambda$, implies that 
\begin{equation*}
\lim_{N\rightarrow \infty}
N\xi^N_{x,x+e_j}(H(\pfrac{x+e_j}{N})-H(\pfrac{x}{N}))
+N\xi^N_{x,x-e_j}(H(\pfrac{x-e_j}{N})-H(\pfrac{x}{N})) = 0 \, .
\end{equation*}
Therefore, the terms in \eqref{inGamma} converges uniformly to zero, and the same holds for the whole sum.
\end{proof}

\smallskip

\begin{corollary}\label{4.4} 
For $H\in \mf D_\Lambda$ and for every $\delta >0$,
 \begin{equation*}
\lim_{N\to\infty} \bb Q_{\mu_N}^{\Lambda,N} \Big[ \sup_{0\le t\le T}
\Big\vert\<\pi^N_t, H\> \,-\, \<\pi^N_0, H \> \,-\,
\int_0^t  \, \<\pi^N_s ,\mc L_{\Lambda}H \> \,ds\,\Big\vert \, > \, \delta \Big] 
\;=\; 0\,.
\end{equation*}
\end{corollary}
\begin{proof}
By a simple calculation, the martingale defined  in \eqref{M1} can be rewritten as
\begin{equation*}
M^{N}_t=\<\pi^N_t, H\> \,-\, \<\pi^N_0, H \> \,-\,
\int_0^t  \,  \<\pi^N_s ,N^2 \bb L_N\,H \> \,ds\,.
\end{equation*}
The result follows from proposition \ref{condicoesI} and expression \eqref{limprob1}.
\end{proof}

\medskip

From Corollary \ref{4.4}, we obtain that (\ref{equation}) in Lemma \ref{c1} holds for any function $H \in C(\bb T^d) \cap \mf D_\Lambda$. This follows from the fact that for such $H$,
$$
\sup_{0\le t\le T} \Big\vert\<\pi^N_t, H\> \,-\, \<\pi^N_0, H \> \,-\,
\int_0^t  \, \<\pi^N_s ,\mc L_{\Lambda}H \> \,ds \Big\vert
$$
is a continuous function in $D([0,T], \mc M)$.

We still need to show that (\ref{equation}) holds for any $H \in \mf D_\Lambda$. For $\varepsilon >0$, define 
$$
(\p\Lambda)^\eps=\{u\in\bb T^d;\,\textrm{dist}(u,\p\Lambda)\leq \eps\} \, .
$$ 
Let $H^\varepsilon$ be a smooth function which coincides with $H$ in $\bb T^2\backslash (\partial\Lambda)^\eps$ and $\sup_{\, \bb T}|H^\varepsilon| \leq \sup_{\, \bb T}|H|$. 

Recall that $\bb Q^*$ is concentrated on trajectories $\pi_t(du) =
\rho(t,u) du$ whose density $\rho$ is non-negative and bounded by 1. Then, under $\bb Q^*$, 
\begin{eqnarray*}
\sup_{0\le t\le T}|\<\pi_t,H^\varepsilon-H\>| & \leq & \sup_{0\le t\le T}\int_{(\partial\Lambda)^\eps}
\!\!\!\rho(t,u)\,|H^\varepsilon(u)-H(u)|\,du \\
& \leq & 2 \, \ell_d((\partial\Lambda)^\eps) \,\sup_{ u \in \bb T}|H(u)|\,.
\end{eqnarray*}
Therefore, for every $\delta >0$,
\begin{equation*}
\bb Q^{*} \Big[ \, \sup_{0\le t\le T}
\Big\vert \<\pi_t, H \> \;-\; 
\<\pi_0,H \> \;-\;
\int_0^t  \,  \<\pi_s ,\mc L_\Lambda H \> \;ds \Big\vert
 \, > \, \delta \Big]
 \end{equation*}
\begin{equation*}
\leq\bb Q^{*} \Big[ \, \sup_{0\le t\le T}
\Big\vert \<\pi_t, H^\varepsilon \> \;-\; 
\<\pi_0,H^\varepsilon \> \;-\;
\int_0^t  \,  \<\pi_s ,\mc L_\Lambda H \> \;ds \Big\vert
 \, > \, \delta/3 \Big]
 \end{equation*}
\begin{equation*}
 +2\,\bb Q^{*} \Big[ \, \sup_{0\le t\le T}
\Big\vert \<\pi_t,H^\varepsilon-H\> \Big\vert
 \, > \, \delta/3 \Big]\,.
\end{equation*}
For small enough $\varepsilon$, the second probability in the sum above is null.

If $G_1$, $G_2$, $G_3$ are continuous functions, the application
from $D([0,T],\mc M)$ to $\bb R$ that associates to a trajectory $\{\pi_t,0\leq t\leq T\}$ the number
\begin{eqnarray*}
 \sup_{0\le t\le T}
\Big\vert \<\pi_t, G_1\> \;-\; \<\pi_0, G_2 \> \;-\;
\int_0^t  \,  \<\pi_s ,G_3 \> \;ds \Big\vert
\end{eqnarray*}
is continuous in the Skorohod metric. Then,

\begin{equation*}
\bb Q^{*} \Big[ \, \sup_{0\le t\le T}
\Big\vert \<\pi_t, H^\varepsilon \> \;-\; 
\<\pi_0,H^\varepsilon \> \;-\;
\int_0^t  \,  \<\pi_s ,\mc L_\Lambda H \> \;ds \Big\vert
 \, > \, \delta/3 \Big]
 \end{equation*}
 \begin{equation*}
 \leq\varliminf_{N\to\infty} \bb Q_{\mu_N}^{\Lambda,N} \Big[ \, \sup_{0\le t\le T}
\Big\vert \<\pi^N_t, H^\varepsilon \> \;-\; 
\<\pi^N_0,H^\varepsilon \> \;-\;
\int_0^t  \,  \<\pi^N_s ,\mc L_\Lambda H \> \;ds \Big\vert
 \, > \, \delta/3\Big]\,,
\end{equation*}
since $ \bb Q_{\mu_N}^{\Lambda,N} $ converges weakly to $ \bb Q^*$ and the above set is open.

By definition,
\begin{equation*}
\sup_{0\le t\le T} \Big\vert \<\pi^N_t,H^\varepsilon-H\> \Big\vert\leq
\frac{1}{N^d}\sum_{x\in\bb T_N^d}\Big\vert H^\varepsilon(x/N)-H(x/N)\Big\vert
\end{equation*}
\begin{equation*}
\leq \Big(\ell_d((\partial\Lambda)^\eps)+O(\pfrac{1}{N})\Big)\,2\,\sup_{ u \in \bb T}|H(u)|\,,
\end{equation*}
because $H^\varepsilon$ coincides with $H$ in $\bb T\backslash (\partial\Lambda)^\eps$. 
Using the same argument as before, we obtain

\begin{equation*}
 \varliminf_{N\to\infty} \bb Q_{\mu_N}^{\Lambda,N} \Big[ \, \sup_{0\le t\le T}
\Big\vert \<\pi_t, H^\varepsilon \> \;-\; 
\<\pi_0,H^\varepsilon \> \;-\;
\int_0^t  \,  \<\pi_s ,\mc L_\Lambda H \> \;ds \Big\vert
 \, > \, \delta/3\Big]
\end{equation*}

\begin{equation*}
 \leq\varliminf_{N\to\infty} \bb Q_{\mu_N}^{\Lambda,N} \Big[ \, \sup_{0\le t\le T}
\Big\vert \<\pi_t, H \> \;-\; 
\<\pi_0,H \> \;-\;
\int_0^t  \,  \<\pi_s ,\mc L_\Lambda H \> \;ds \Big\vert
 \, > \, \delta/9\Big]
\end{equation*}

\begin{equation*}
 +2\varliminf_{N\to\infty} \bb Q_{\mu_N}^{\Lambda,N} \Big[ \, \sup_{0\le t\le T}
\Big\vert \<\pi_t,H^\varepsilon-H\> \Big\vert
 \, > \, \delta/9 \Big]\,.
\end{equation*}
Again, for small enough $\varepsilon$, the second probability in the sum above is null.
Recalling the corollary \ref{4.4},  we finally conclude that $\bb Q^{*}$ is concentrated on absolutely continuous paths
$\pi_t(du) = \rho(t,u) du$ with positive density bounded by $1$, and $\bb Q^{*}$ a.s. 
\begin{equation*}
\dl \rho_t,  H \dr - \dl \rho_0, H \dr \;=\; \int_0^t  \, \dl \rho_s \,,\, \mc L_\Lambda  H \dr \;ds\,,
\end{equation*}
\noindent for any $H\in\mf D_\Lambda$. Therefore we have proved Lemma \ref{c1}.


\subsection{Uniqueness of weak solutions}\label{uniqueness}
Now, we prove that the solution of (\ref{edp}) is unique.
It suffices to check that the only solution of (\ref{edp}) with $\rho_0\equiv 0$ is $\rho\equiv 0$, 
because of the linearity of $\mc L_\Lambda$. Let 
 $\rho :\bb R_+ \times \bb T^d \to \bb R$ be a weak solution of the parabolic differential equation
\begin{equation*}
\left\{
\begin{array}{l}
{\displaystyle \partial_t \rho \; =\; \mc L_\Lambda \rho } \\
{\displaystyle \rho(0,\cdot) \;=\; 0\,.}
\end{array}
\right.
\end{equation*} By definition,
\begin{equation}\label{eqint}
\dl \rho_t, H\dr
\;=\; \int_0^t \dl \rho_s , \mc L_\Lambda  H \dr \, ds\;,
\end{equation}
 for all functions $H$ in $\mc H^1_\Lambda $ and all $t>0$. 
 From the theorem \ref{Lambda}, 
the operator  $- \mc L_\Lambda$  has countable eigenvalues  $\{\mu_n : n\ge 0\}$ and  eigenvectors $\{F_n\}$.
 All eigenvalues have finite multiplicity,
  $0= \mu_0 \le \mu_1 \le \cdots$, and $\lim_{n\to\infty} \mu_n  = \infty$.
Besides, the eigenvectors $\{F_n\}$ form a complete orthonormal system in the $L^2(\bb T^d)$. Define 
\[R(t)=\sum_{n\in\bb{N}}\frac{1}{n^{2}(1+\mu_n)}\dl\rho_t,F_n\dr^{2},\] for all $t>0$. Notice that $R(0)=0$
and  $R(t)$ is well defined because $\rho_t$ belongs to $L^{2}(\bb T^d)$. 
Since $\rho$ satisfy \eqref{eqint}, we have that $\frac{d}{dt}\dl\rho_t,F_n\dr^{2}=
-2\mu_n\dl\rho_t,F_n\dr^{2}$. Then
\begin{equation*}
(\pfrac{d}{dt}R)(t)=-\sum_{n\in\bb N}\frac{2\mu_n}{n^{2}(1+\mu_n)}\dl \rho_t,F_n\dr ^{2}\,,
\end{equation*}
because 
$\sum_{n\leq N}\frac{-2\mu_n}{n^{2}(1+\mu_n)}\dl\rho_t,F_n\dr^{2}$ converges uniformly to
$\sum_{n\in \bb N}\frac{-2\mu_n}{n^{2}(1+\mu_n)}\dl\rho_t,F_n\dr^{2}$, 
when $N$ increases to infinity. Thus $R(t)\geq0$  and $(\frac{d}{dt}R)(t)\leq 0$, for all $t>0$ and $R(0)=0$. 
From this, we obtain $R(t)=0$ for all $t>0$. 
Since $\{F_n\}$ is a complete orthonormal system, $\<\rho_t,\rho_t\>=0$, for all $t>0$, 
which implies $\rho\equiv0$.

\bigskip

\noindent {\emph{Acknowledgments:}} The authors T. Franco and A. Neumann would like to thank
Claudio Landim, their PhD advisor, for support and valuable comments.


\medskip

\begin{thebibliography}{99}

\bibitem{ct} A. Chambolle, G. Thouroude: {\em Homogenization of interfacial energies and construction of plane-like minimizers in periodic media through a cell problem}, Netw. Heterog. Media  {\bf 4}, $n^o$ 1 127--152 (2009).

\bibitem{e} L. \ Evans, {\em Partial Differential Equations}.  [Graduate Studies in Mathematics], American Mathenatical Society.

\bibitem{f} A. Faggionato. {\em Bulk diffusion of 1D exclusion process with bond disorder}. Markov Processes and
Related Fields 13, 519-542 (2007).

\bibitem{f2} A. Faggionato. {\em Hydrodynamic limit of symmetric exclusion processes
in inhomogeneous media. ArXiv, http://arxiv.org/pdf/1003.5521v1} (2010).

\bibitem{fjl} A. Faggionato, M. Jara, C. Landim: {\em Hydrodynamic behavior of one dimensional subdiffusive exclusion processes with random conductances}, Probab. Th. and Rel. Fields {\bf 144}, $n^o$ 3-4, 633--667 (2008).

\bibitem{fl} T. Franco, C. Landim: {\em Hydrodynamic Limit of Gradient Exclusion Processes with Conductances}. Archive for Rational Mechanics and Analysis (Print), v. 195, p. 409-439, (2010). 

\bibitem{j} M. Jara, {\em Hydrodynamic limit of particle systems in inhomogeneous media}, online,  ArXiv	 http://arxiv.org/abs/0908.4120

\bibitem{kl} C.\ Kipnis, C.\ Landim, {\em Scaling limits of interacting
  particle systems}. Grundlehren der Mathematischen Wissenschaften
  [Fundamental Principles of Mathematical Sciences], 320.
  Springer-Verlag, Berlin (1999).

\bibitem{n} K.\ Nagy, {\em Symmetric random walk in random
    environment}. Period. Math. Ung. {\bf 45}, 101--120 (2002).

\bibitem{v} F. Valentim: {\em Hydrodynamic limit of gradient exclusion processes with conductances on $\mathbb{Z}^d$}, preprint.

\bibitem{z} E. Zeidler, {\em Applied Functional Analysis. Applications
    to Mathematical Physics.}. Applied Mathematical Sciences, 108.
  Springer-Verlag, New York (1995).

\end{thebibliography}
\end{document}